\documentclass[11pt]{article}

\usepackage{amsthm}
\usepackage{fullpage}

\newtheorem{theorem}{Theorem}

\newtheorem{corollary}{Corollary}

\newtheorem{definition}{Definition}

\newtheorem{lemma}{Lemma}

 \RequirePackage{ifthen}
\RequirePackage{times}

\RequirePackage{pifont}

\usepackage{amssymb}
\usepackage{amsmath}
\usepackage{graphicx}

\usepackage{epsfig}

\newcommand\Li{\mathrm{Li}}

\begin{document}

\title{
 Polynomials associated with
Partitions: \
Their Asymptotics and Zeros
}

\author{Robert  P. Boyer, William M. Y. Goh}

\maketitle

\begin{abstract}
Let $p_n$ be the number of partitions of an integer $n$.
For each of the partition statistics of counting their parts, ranks, or cranks,
there is a natural family of integer polynomials.
We  investigate  their 
asymptotics and the limiting behavior of their zeros as sets and densities.
   \end{abstract}

 \section{Introduction}\label{section:introduction}

The purpose of this paper is to survey
several natural polynomial families associated with integer partitions
focusing on
 their  asymptotics and the limiting behavior
of their zeros. Our principal families are
 
\begin{enumerate}
\item
Taylor polynomials of the analytic function 
$P(x)=\prod_{n\geq 1} (1-x^n)^{-1}$,
the generating function of the partition numbers
(Section \ref{section:Taylor}).
\item
Polynomials $F_n(x)$ associated with counting partitions in parts
(Section \ref{section:stanley}).
\item
Polynomials associated with the rank or crank of a partition
(Section \ref{section:crank}).
\end{enumerate}

We introduce several definitions used throughout the paper.  

\begin{definition}\label{def:zero_attractor}
Let ${\mathcal Z}(q_n)$ denote the finite set of zeros of
the polynomial $q_n$. Then the {\sl zero attractor} $\mathcal A$
of the polynomial sequence $\{q_n\}$ whose degrees go to $\infty$
 is the limit of
${\mathcal Z}(q_n)$ in the Hausdorff metric $\Delta$
on the non-empty compact
subsets $\mathcal K$ of ${\mathbb C}$.
\end{definition}

We recall the standard:

\begin{definition}\label{def:asymptotic_zero}
The {\sl asymptotic zero distribution} 
for a sequence $\{q_n\}$ of polynomials 
whose degrees go to $\infty$
is the weak$^{*}$-limit
of the normalized counting measures of their zeros
$\frac{1}{\deg(q_n)} \sum \{ \delta_z : z \in {\mathcal Z}(q_n)\}$.
\end{definition}

We single out a useful compromise from the obtaining the
full asymptotic zero distribution.

\begin{definition}\label{def:argument}
We say that the  arguments of the zeros of a polynomial family $\{q_n(x)\}$ 
whose degrees go to $\infty$
are 
{\sl uniformly distributed on
the unit circle as $n\to\infty$} if the normalized counting measures
$\frac{1}{\deg(n)} \sum \{ \delta_{\arg z} : z \in {\mathcal Z}(q_n)\}$ converge in the
weak$^{*}$-topology to normalized Lebesgure measure on the unit circle.
\end{definition}

The following   result of Erd\"os and Tur\'an 
(\cite{erdos}, Theorem 1)
will be used repeatedly throughout the paper to
determine that the arguments of zeros are uniformly distributed.
Let $q(x)$ be the polynomial
$\sum_{k=0}^n a_k x^k$ of degree $n$ 
 with non-zero constant term $a_0 \neq 0$.
 For $0 \leq \theta_1 < \theta_2 \leq 2 \pi$,
 \begin{equation}\label{eq:erdos}
 \left|
 \# \, \{ z : \arg z \in [\theta_1,\theta_2], q(z) =0 \}
 - \frac{ \theta_2-\theta_1}{2 \pi} n
 \right| < 
 16 \sqrt{
 n \ln \left(
\frac{ |a_0| + |a_1| + \cdots + |a_n|} { \sqrt{ a_0 a_n}}
 \right)
 }.
 \end{equation}

\section{Taylor Polynomials of $P(x)$}\label{section:Taylor}

Let $p_k$ be the number of partitions of a positive integer $k$ with $p_0=1$
by convention. The ordinary generating function $P(x)$ for $\{p_k\}$ is
\begin{equation}
P(x)=
\prod_{n\geq 1} \frac{1}{1-x^n} = \sum_{k=0}^\infty p_k x^k.
\end{equation}
A natural choice of polynomials associated with the partitions is simply
the Taylor polynomials $s_n(x)$ of $P(x)$:
\begin{equation}
s_n(x) = \sum_{k=0}^n p_k x^k
\end{equation}
since $P(x)$ is analytic in the open unit disk $\mathbb D$.

The asymptotics of these polynomials $s_n(x)$ depend 
on the classical result of the asymptotics of the partition
numbers $p_n$:
\begin{equation}\label{eq:hardy_ramanujan}
 p_n \left/  \left[   \frac{1}{4n\sqrt{3}} 
\exp \left( \pi  \sqrt{\frac{2n}{3}}  \right) \right]  \right. \to 1.
\end{equation}
See either \cite{andrews} or
 \cite{ayoub}.
We first establish the limiting behavior of their zeros.

\begin{theorem}
(a) The zero attractor of the Taylor polynomials 
$\{ s_n(x)\}$ is the unit circle.
\\
(b) The asymptotic zero density is Lebesgue measure on the unit
circle.
\end{theorem}

\begin{proof}
Recall the Enestr\"om-Kakeya Theorem:
If the coefficients of the polynomial
 $q(z)= \sum_{k=0}^n a_k z^k$ 
satisfy $a_n \geq a_{n-1} \geq \cdots \geq a_0 \geq 0$, then all
the zeros of $p(z)$ lie in the closed unit disk (see \cite{marden}, p. 136).
Since the partition numbers are positive and increasing, 
the zeros of the Taylor polynomials $s_n(x)$ must lie in 
the closed unit disk
$\overline{{\mathbb D}}$. 

Next let $f(x) = \sum_{k=0}^\infty c_k x^k$ be an analytic
function with radius of convergence 1.
To state the Jentzsch Theorem \cite{erdos}
concerning the
zeros of the Taylor polynomials $t_n(x)$ of $f(x)$, 
recall that  $a$ is called a
 limit point of zeros of
$t_n(x)$ if for every $\varepsilon >0$ there are infinitely many indices
$n$ so $t_n(z_n)=0$ with $|z_n- a|< \varepsilon$.  
Then the collection of all
limit points of zeros of $t_n(x)$ must contain the unit
circle.

Since $s_n(x)$ are the Taylor polynomials
of the generating function $P(x)$ 
which is analytic and
does not vanish in $\mathbb D$, no limit point of 
the polynomials $s_n(x)$ can
lie inside $\mathbb D$ since such a limit point must be a zero of $P(x)$.
Since the radius of convergence of $P(x)$ is 1, we conclude
 that  the limit points are exactly the unit circle.
 We conclude that the zero attractor is the unit circle since the
 all the zeros of $s_n(x)$ are bounded in modulus by 1.

  For the polynomials $s_n(x)$, their constant terms are always 1 while their coefficients
 are all bounded above by $p_n$, so the right-hand side of inequality of
 Erd\"os-Tur\'an  (\ref{eq:erdos}) is dominated by
 $16 \sqrt{ n \ln( n \sqrt{ p_n})}$.
 Hence, we have the following limit by (\ref{eq:hardy_ramanujan}):
 \begin{eqnarray*}
\left|
 \frac{1}{n} \# \{ z : \arg z \in [\theta_1,\theta_2], q(z) =0 \}
 - \frac{ \theta_2-\theta_1}{2 \pi} 
 \right|<
16 \sqrt{ \frac{1}{n} \ln( n \sqrt{ p_n})} \to 0.
 \end{eqnarray*}
 A compactness argument shows that the unit circle is the zero attractor.
\end{proof}

Because of the non-negativity and monotonicity of the coefficients of 
$s_n(x)$ together with
 the subexponential growth 
 of $s_n(1)$, both the zero attractor and the asymptotic
zero distribution for $s_n(x)$
were quickly obtained. A more complete understanding
of these polynomials, though, requires 
their asymptotics outside the unit disk.
In general, it is very useful to have asymptotic expansions for a polynomial
family throughout the complex plane.  In \cite{boyer_goh_euler}, we obtained
such expansions for the Euler and Bernoulli polynomials 
while in Section
\ref{section:stanley} we describe expansions for another partition polynomial
family. Further,  we note that the Euler and Bernoulli polynomial zeros are not
uniformly distributed around the unit circle and that the zero distribution studied
in Section \ref{section:stanley} is more subtle than any of these examples.

\begin{theorem}
Let $\delta >0$  and 
$0 < \eta < 1/2$,
 then
\begin{equation*}
s_n(x)=\frac{x^{n+1}}{x-1}\frac{e^{a\lambda _{n}}
\,
\lambda _{n}^{-2}}{4\sqrt{3}}
\,
\left(1+O_{\delta }(\lambda _{n}^{ -\eta }) \right),
\end{equation*}
where 
\begin{equation}
a=\pi \sqrt{ 2/3 }, \quad   \lambda _{n}= \sqrt{n - 1/24}
 \label{lamdan}
\end{equation}
 and the
constant in the big oh term $O_{\delta }(\lambda _{n}^{-\eta })$
 depends only on $\delta $  and holds uniformly for all $x$
with $|x| \geq 1+\delta$.
\end{theorem}

\begin{proof} 
For any $0<r<1$, we have
\begin{equation*}
p_{n}=\frac{1}{2\pi i}
\oint_{\left| \zeta \right| =r}\frac{P(\zeta )}{\zeta^{n+1}}  \,d\zeta .
\end{equation*}
By summing over the above expression for the partition numbers
$p_n$, we obtain an integral form for the
Taylor polynomial:
\begin{eqnarray*}
s_n(x)=\frac{1}{2\pi i}
\oint_{\left| \zeta \right| =r}
\frac{P(\zeta )}{\zeta }
\left(
\sum_{j=0}^{n}(\frac{x}{\zeta })^{j}
\right)\, d\zeta 
=
\frac{1}{2\pi i}
\oint_{\left| \zeta \right| =r}
\frac{P(\zeta )}{\zeta -x}
\left(1-   \left( \frac{x}{\zeta } \right)^{n+1} \right)  \,d\zeta .
\end{eqnarray*}
Since $\left| x\right| \geq 1+\delta$,
we find, by using the Cauchy integral theorem, that 
\begin{equation*}
\frac{1}{2\pi i}
\oint_{\left| \zeta \right| =r}
\frac{P(\zeta )}{\zeta -x}
\,d\zeta =0.
\end{equation*}
The integral for the Taylor polynomial $s_n(x)$ reduces to
\begin{equation}
s_n(x)=
\frac{-x^{n+1}}{2\pi i}
\oint_{\left| \zeta \right| =r}
\frac{P(\zeta)}{\zeta -x}\zeta^{-n-1}\,d\zeta.  \label{vn}
\end{equation}
Next we define $I_n$ as the integral:
\begin{equation}
I_{n}=\frac{1}{2\pi i}
\oint_{\left| \zeta \right| =r}
\frac{P(\zeta )}{\zeta-x}\zeta^{-n-1}\,d\zeta.   \label{in}
\end{equation}
In particular, we must have:
\begin{equation}
s_n(x)=-x^{n+1}I_{n}.  \label{vnin}
\end{equation}

Thus our goal is to find an asymptotic approximation for $I_{n}$.
 Our
strategy follows very closely  that of \cite{ayoub}. 
Consequently, we will adopt
the same notation as  Ayoub to avoid confusion. 
Not surprisingly, the methods come from a proof of the asymptotics
of partition numbers originally by J.  Upsensky and uses the functional
equation of the modular function.
Basically, the major
contribution to the integral in (\ref{in}) comes from the a small
neighborhood of the strongest singularity $\zeta =1$ of $P(\zeta )$.

We begin with  the following well-known functional equation
which  is essential: 
For $\Re(\tau )>0$,
\begin{equation}
P(e^{-2\pi \tau })
=
\psi (\tau )P(e^{-2\pi /\tau }),  \label{functional}
\end{equation}
where
\begin{equation*}
\psi (\tau )=\sqrt{\tau }
\exp \left[\frac{\pi }{12} \left(\frac{1}{\tau }-\tau \right) \right].
\end{equation*}
Now we put
$\zeta =e^{-2\pi \tau }$ with $ d\zeta =  e^{-2\pi \tau }  2\pi i  \, d\phi$
where 
$\tau =\alpha -i\phi$,
with $\alpha =\alpha (n)>0$.
 Note that we shall choose $\alpha $ so that 
$\alpha \rightarrow 0$ as $n\rightarrow \infty$.
The specific form of $\alpha $ will be made clear below.
Using the functional equation, we write $I_n$ as:
\begin{eqnarray}
I_{n} =
\int_{-1/2}^{1/2}
\frac{P(e^{-2\pi \tau })}{e^{-2\pi \tau }-x}
e^{2\pi n\tau } \,d\phi 
=
J+\widetilde{I_{n}},  \label{inin}
\end{eqnarray}
where 
\begin{equation*}
J=\int_{-1/2}^{1/2}\frac{\psi (\tau )}{e^{-2\pi \tau }-x}e^{2\pi n\tau } \, d\phi,
\quad
\widetilde{I_{n}}=\int_{-1/2}^{1/2}
\frac{P(e^{-2\pi \tau })-\psi (\tau )}{e^{-2\pi \tau }-x}e^{2\pi n\tau } \, d\phi .
\end{equation*}

To estimate $\widetilde{I_{n}}$, we break the interval into three parts; a
neighborhood of origin, say $-\phi _{0}\leq \phi \leq \phi _{0},$ and the
remaining two segments from $-1/2$ to $-\phi _{0}$ and $\phi _{0}$ to $1/2$.
Choose 
$\phi _{0}=\lambda \alpha $
and $\lambda$  that satisfies 
$2\pi =\alpha (1+\lambda ^{2})$,
  that is,  
$\phi _{0}=(2\pi \alpha -\alpha^{2})^{1/2}$.

We proceed as in \cite{ayoub} to get the estimates:

\begin{lemma}
(a) For $\left| \phi \right| \leq \phi _{0}$ we have
\begin{equation}
P(e^{-2\pi \tau })-\psi (\tau )=O(1).  \label{claima}
\end{equation}

(b) For $\phi _{0}\leq \phi \leq \frac{1}{2}$ or 
$\frac{-1}{2}\leq \phi \leq
-\phi _{0}$ we have
\begin{equation}
P(e^{-2\pi \tau })-\psi (\tau )=O(e^{\pi /(48\alpha )}).  \label{claimb}
\end{equation}
\end{lemma}

\begin{proof}
 Equation (\ref{functional}) is essential for the proof here. For
details see  equations (14) and (19) on page 150 of \cite{ayoub}. 
\end{proof}

We use this lemma to estimate $\widetilde{I_{n}}$.
Define $\widetilde{I}_{n,1}$, $\widetilde{I}_{n,2}$, and
$\widetilde{I}_{n,3}$ as:
\begin{eqnarray*}
\widetilde{I_{n}}
&=&
\left(\int_{-1/2}^{-\phi _{0}}+
\int_{-\phi _{0}}^{\phi_{0}}+\int_{\phi _{0}}^{1/2}
\right)
\,
\left(
\frac{P(e^{-2\pi \tau })-
\psi (\tau )}{e^{-2\pi \tau }-x}
e^{2\pi n\tau }
\right)
 \,d\phi 
\\
&=&
\widetilde{I}_{n,1}+\widetilde{I}_{n,2}+\widetilde{I}_{n,3}.
\end{eqnarray*}

From equation (\ref{claima}),
\begin{equation*}
\widetilde{I}_{n,2}
=O
\left(
\int_{-\phi _{0}}^{\phi _{0}}
\left| \frac{1}{e^{-2\pi \tau }-x}\right| \left| e^{2\pi n\tau }\right| \, d\phi 
\right).
\end{equation*}
For $\left| x\right| \geq 1+\delta $
\begin{equation*}
\left| \frac{1}{e^{-2\pi \tau }-x}\right| 
\leq 
\frac{1}{\left| x\right|
-\left| e^{-2\pi \tau }\right| }
=
\frac{1}{\left| x\right| -e^{-2\pi \alpha }}
\leq 
\frac{1}{\left| x\right| -1}
\leq
 \frac{1}{\delta }.
\end{equation*}

Hence
$
\widetilde{I}_{n,2}
=O_{\delta }(e^{2\pi n\alpha }),
$
whereas from equation (\ref{claimb}) 
$
\widetilde{I}_{n,3}=O_{\delta } \left(e^{2\pi n\alpha +\pi /(48\alpha )} \right);
$
and exactly the same estimate holds for $\widetilde{I}_{n,1}$.

 From  equation (\ref{inin}), we have now shown the following:
 
\begin{lemma}
\begin{equation}
I_{n}=J+O_{\delta }\left(e^{2\pi n\alpha +\pi /(48\alpha )} \right),  \label{inj}
\end{equation}
\end{lemma}

We use the functional equation
(\ref{functional}) to obtain
\begin{equation*}
J=\int_{-1/2}^{1/2}\frac{(\alpha -i\phi )^{1/2}}{e^{-2\pi \tau }-x}
\exp \left(\frac{\pi }{12(\alpha -i\phi )}+2\pi (n-\frac{1}{24})(\alpha -i\phi ) \right)
\,  d\phi .
\end{equation*}
For convenience, we put 
\begin{equation}
m=2\pi \left(n-\frac{1}{24} \right)
=2\pi \lambda _{n}^{2}.  \label{m}
\end{equation}


We change  variables $\phi =\alpha u$ to get 
\begin{equation*}
J
=
\alpha^{3/2}\int_{-1/(2\alpha )}^{1/(2\alpha )} \,
\frac{(1-iu)^{1/2}}{e^{-2\pi \alpha (1-iu)}-x} \,
\exp \left(
\frac{\pi }{12\alpha (1-iu)}
+
m\alpha
(1-iu)
\right)  \, du.
\end{equation*}
To obtain an asymptotic approximation for $J,$ we set the coefficients of 
$\frac{1}{1-iu}$ and $1-iu$ to be equal. Thus 
$
\frac{\pi }{12\alpha }=m\alpha =\sigma,
$
and so
\begin{equation}
\alpha =\sqrt{\frac{\pi }{12m}},
\quad 
\sigma =\sqrt{\frac{\pi m}{12}},
\label{alpsig}
\end{equation}
where $m$ was defined in equation (\ref{m}). This is how $\alpha $ is made
explicit. 
Consequently, 
\begin{eqnarray}
J
&=&
\alpha ^{3/2}\int_{-1/(2\alpha )}^{1/(2\alpha )}
\frac{(1-iu)^{1/2}}{e^{-2\pi \alpha (1-iu)}-x}
\exp \left[
\, \sigma \left(\frac{1}{1-iu}+(1-iu)  \right) 
\right] \,du.
\nonumber
\\
&=&
\alpha^{3/2}e^{2\sigma }\int_{-u_{0}}^{u_{0}}
\frac{(1-iu)^{1/2}}{e^{-2\pi\alpha (1-iu)}-x}\exp [-\sigma g(u)] \, du,  \label{j}
\end{eqnarray}
where 
\begin{equation*}
g(u)=\frac{u^{2}}{1-iu}, \quad 
u_{0}=\frac{1}{2\alpha }.
\end{equation*}
Note that from equation (\ref{alpsig}) 
$\sigma =\sigma (n)\rightarrow \infty $ 
and $\alpha \rightarrow 0$ as $n\rightarrow \infty$.

 To approximate $J$
we follow   (\cite{copson},page 91). We choose $\varepsilon$
to lie in the interval $(1/3,1/2)$.
Write
\begin{eqnarray}
J&=&
\alpha ^{3/2}e^{2\sigma }
\left[
\int_{-u_{0}}^{-\sigma^{-\varepsilon}}
+\int_{-\sigma^{-\varepsilon }}^{\sigma^{-\varepsilon }}
+
\int_{\sigma^{-\varepsilon }}^{u_{0}} \right]
\frac{(1-iu)^{1/2}\exp [-\sigma g(u)]}
{e^{-2\pi\alpha (1-iu)}-x}  \, du
\nonumber 
\\
&=&
J_{1}+J_{2}+J_{3}.  \label{j1j2j3}
\end{eqnarray}

\begin{lemma}
(a)
\qquad
Both $J_1$ and $J_3$ equal
$\displaystyle \frac{ \alpha^{3/2} e^{2 \sigma} \sqrt{\pi}}{ \sigma} o_\delta( \sigma^{1-3\epsilon})$.
\\
(b)
\qquad
$\displaystyle J_2 
=
\frac{ \alpha^{3/2} e^{2 \sigma}}{1-x} \, \frac{ \sqrt \pi}{ \sqrt \sigma}
( 1+ O_\delta( \sigma^{1-3 \epsilon}))$.
\end{lemma}

\begin{proof}
We estimate $J_{2}$ first. Note for 
$-\sigma^{-\varepsilon }\leq u\leq
\sigma^{-\varepsilon }$ we have

\begin{eqnarray*}
(1-iu)^{1/2}
&=&
1+O(\sigma ^{-\varepsilon })=1+O(\sigma ^{1-3\varepsilon })
\\
\frac{1}{e^{-2\pi \alpha (1-iu)}-x}
&=&
\frac{1}{1-x+O(\alpha )}
=
\frac{1}{1-x}
\frac{1}{1+O_{\delta }(\alpha )}
\\
&=&
\frac{1}{1-x}(1+O_{\delta }(\alpha ))
=\frac{1}{1-x}(1+O_{\delta }(\sigma^{-1}))=\frac{1}{1-x}
(1+O_{\delta }(\sigma^{1-3\varepsilon })),
\\
g(u)
&=&
\frac{u^{2}}{1-iu}=u^{2}+O(\sigma^{-3\varepsilon })
\end{eqnarray*}
so that 
\begin{equation*}
\exp [-\sigma g(u)]=\exp [-\sigma u^{2}](1+O(\sigma^{1-3\varepsilon })).
\end{equation*}
Making the above substitutions, we find 
\begin{eqnarray}
J_{2}
&=&
\alpha ^{3/2}e^{2\sigma }
\int_{-\sigma ^{-\varepsilon }}^{\sigma^{-\varepsilon }}
\frac{(1-iu)^{1/2}\exp [-\sigma g(u)]}{e^{-2\pi \alpha
(1-iu)}-x} \, du
\\
&=&
\alpha ^{3/2}e^{2\sigma }
\int_{-\sigma^{-\varepsilon }}^{\sigma^{-\varepsilon }}
\frac{\exp [-\sigma u^{2}]}{1-x}(1+O_{\delta }
(\sigma^{1-3\varepsilon }))\,du
\\
&=&
\frac{\alpha^{3/2}e^{2\sigma }}{1-x}
\left(
\int_{-\sigma^{-\varepsilon}}^{\sigma^{-\varepsilon }}e^{-\sigma u^{2}}  \, du
\right)
\,
\left(
1+O_{\delta }(\sigma^{1-3\varepsilon }) \right)  \label{j2}
\end{eqnarray}
Now 
\begin{equation*}
\int_{-\sigma ^{-\varepsilon }}^{\sigma ^{-\varepsilon }}
e^{-\sigma u^{2}}  \,du
=
\int_{-\infty }^{\infty }e^{-\sigma u^{2}} \,du
-
\left[
\int_{-\infty}^{-\sigma^{-\varepsilon }}
+
\int_{\sigma^{-\varepsilon }}^{\infty}
\right]
e^{-\sigma u^{2}}  \, du.
\end{equation*}
It is not hard to see that 
since
$\int_{-\infty }^{\infty } e^{-\sigma u^{2}}  \, du
=
{\sqrt{\pi }}/{\sqrt{\sigma }} 
$
\begin{equation*}
\left[
\int_{-\infty }^{-\sigma^{-\varepsilon }}
+
\int_{\sigma^{-\varepsilon}}^{\infty } \right]
\,
e^{-\sigma u^{2}} \,du
=o(\sigma^{1-3\varepsilon })
\end{equation*}
so that 
\begin{equation*}
\int_{-\sigma^{-\varepsilon }}^{\sigma ^{-\varepsilon }}
e^{-\sigma u^{2}}  \,du=
\frac{\sqrt{\pi }}{\sqrt{\sigma }}(1+o(\sigma^{1-3\varepsilon })).
\end{equation*}
Hence from equation (\ref{j2}) we get 
\begin{eqnarray*}
J_{2}
&=&
\frac{\alpha ^{3/2}e^{2\sigma }}{1-x}\frac{\sqrt{\pi }}{\sqrt{\sigma }}
\left(1+o(\sigma ^{1-3\varepsilon }) \right)
\,
\left(1+O_{\delta }(\sigma ^{1-3\varepsilon }) \right)
\\
&=&
\frac{\alpha^{3/2}e^{2\sigma }}{1-x}\frac{\sqrt{\pi }}{\sqrt{\sigma }}
\left(
1+O_{\delta }(\sigma ^{1-3\varepsilon })
\right).
\end{eqnarray*}
Recall 
\begin{equation*}
J_{3}=\alpha^{3/2}e^{2\sigma }
\int_{\sigma^{-\varepsilon }}^{u_{0}}
\frac{(1-iu)^{1/2}\exp [-\sigma g(u)]}{e^{-2\pi \alpha (1-iu)}-x}  \, du.
\end{equation*}
We have the estimates
\begin{eqnarray}
\left| J_{3}\right| 
&\leq&
 \alpha^{3/2}e^{2\sigma }
\int_{\sigma^{-\varepsilon }}^{u_{0}}
\left|
 \frac{(1-iu)^{1/2}
\exp \left[-\sigma g(u) \right]}{
e^{-2\pi \alpha (1-iu)}-x}
\right| \,  du
\nonumber
\\
&=&
\alpha^{3/2}e^{2\sigma }
\int_{\sigma ^{-\varepsilon }}^{u_{0}}
\frac{
(1+u^{2})^{1/4}
\exp \left[-\sigma \, \Re  g(u) \right]}
{\left| e^{-2\pi \alpha (1-iu)}-x\right| }  \,du
\nonumber
\\
&\leq&
 \frac{\alpha^{3/2}e^{2\sigma }}{\delta }
\int_{\sigma^{-\varepsilon}}^{u_{0}}(1+u^{2})^{1/4}
\exp \left[-\sigma \frac{u^{2}}{1+u^{2}} \right] \,du.  \label{j3}
\end{eqnarray}
Since  ${u^{2}}/{(1+u^{2})}$ is an increasing function of $u$,
 we have,
for $u_{0}\geq u\geq \sigma^{-\varepsilon }$,
$
{u^{2}}/{(1+u^{2})}
\geq
{\sigma^{-2\varepsilon }}/{(1+\sigma^{-2\varepsilon })}.
$
This implies 
$
\exp \left( -\sigma {u^{2}}/{(1+u^{2})}  \right)
\leq 
$
$
\exp \left(  -{\sigma^{1-2\varepsilon }}/
{(1+\sigma^{-2\varepsilon })}  \right)
$
$
\leq
 \exp \left(  -{\sigma^{1-2\varepsilon }}/{2}  \right).
$
By assumption $1/3< \varepsilon <1/2$, so we find that
$\exp (-\frac{\sigma^{1-2\varepsilon }}{2})$ is
much smaller than $\sigma ^{1-3\varepsilon }$. 
Hence by the inequality (\ref{j3}) we get
\begin{equation*}
J_{3}
=
\frac{\alpha^{3/2}e^{2\sigma }
\sqrt{\pi }}{\sqrt{\sigma }}o_{\delta}(\sigma^{1-3\varepsilon }).
\end{equation*}
Exactly the same estimate holds for $J_{1}$.
\end{proof}

We now return to the proof of the Theorem.
By the definition of $J_1$, $J_2$, $J_3$ (see equation (\ref{j1j2j3})),
we see that
\begin{equation*}
J=\frac{\alpha^{3/2}e^{2\sigma }}{1-x}\frac{\sqrt{\pi }}{\sqrt{\sigma }}
 \left(1+O_{\delta }(\sigma ^{1-3\varepsilon })  \right).
\end{equation*}
From equation (\ref{inj}) 
\begin{equation*}
I_{n}
=
\frac{\alpha^{3/2}e^{2\sigma }}{1-x}\frac{\sqrt{\pi }}{\sqrt{\sigma }}
(1+O_{\delta }(\sigma^{1-3\varepsilon }))+O_{\delta }(e^{2\pi n\alpha 
+
\pi/(48\alpha )}).
\end{equation*}
To see the final result, we recall the equations
 (\ref{m}), (\ref{lamdan}), and (\ref{alpsig}). It is
convenient that we express everything in terms of $\lambda_{n}$ which 
equals  $\sqrt{n-  1/24}$.
 Thus, with $a=\pi \sqrt{2/3}$, 
\begin{eqnarray*}
\alpha &=&\frac{1}{2\sqrt{6}\lambda _{n}}, \quad
\sigma 
=
\frac{\pi \lambda _{n}}{\sqrt{6}},
\\
2\sigma 
&=&
a\lambda _{n}.
\\
2\pi n\alpha +\pi /(48\alpha )
&=&
\frac{\pi n}{\sqrt{6}\lambda _{n}}+\frac{\pi
\lambda _{n}}{4\sqrt{6}}
\\
&=&
\frac{5\pi \lambda _{n}}{4\sqrt{6}}+o(1)=\frac{5a\lambda _{n}}{8}+o(1).
\end{eqnarray*}
Since $e^{{5a\lambda _{n}}/{8}}$ is dominated by 
$e^{a\lambda _{n}},$ we
have
\begin{equation*}
I_{n}
=
\frac{e^{a\lambda _{n}}\lambda _{n}^{-2}}{(1-x)4\sqrt{3}}
\left(1+O_{\delta}(\lambda _{n}^{1-3\varepsilon })  \right).
\end{equation*}
By comparing with equation (\ref{vnin}) and setting
$\eta = 1- 3 \varepsilon$, we find that 
 the proof is complete.
\end{proof} 


\section{Polynomials for Partitions with Parts}\label{section:stanley}

Let $p_k(n)$ denote the number of partitions of $n$ with exactly $k$ parts.
Define the polynomials $F_n(x) = \sum_{k=1}^n p_k(n) x^k$, the
{\sl partition with parts polynomials}.
They have generating function:
\[
 P(x,u)=\prod_{k \geq 1} \frac{1}{1-xu^k}
=
\sum_{n=1}^\infty F_n(x) u^n.
\]
With $x=1$, $P(1,u)$ reduces to the generating function $P(x)$ for the
partition numbers.
To calculate these polynomials, we make use of the recurrence
$
p_k(n)=p_{k-1}(n-1) + p_{k}(n-k).
$
and the fact that  about half their coefficients  are actually given 
by the partition numbers: $p_{n-k}(n) = p(k)$, $2k <n-1$.

It is also known that the coefficients of $F_n(x)$ are unimodal
for $n$ sufficiently large (\cite{andrews}, page 100).
These polynomials are mentioned in \cite{durfee_poly} where it is pointed out
that they have complex zeros.
Unfortunately, these facts do not  give a hint to the complexity of
their zeros (see Figure 2b).
In fact, Richard Stanley plotted the zeros of $F_{200}(x)$ and asked
what happens at $n\to \infty$.
The proofs of the following results are found in \cite{boyer_goh}.

The asymptotics for $F_n(x)$ outside the unit disk can be found using
the method of Darboux.  We state:

\begin{theorem}
\label{fnx-out}
On compact subsets $K$ that lie in the open set $\{ z : |z| >1\}$,
the polynomials $F_n(x)$ have the asymptotic form
\begin{equation*}
F_{n}(x)=
x^{n}P\left(1,\frac{1}{x} \right)
+O(\left| x\right|^{  C  n} ),  \label{fnx-out_2}
\end{equation*}
where $1/2<C <1$ and
the big $O$ term holds uniformly in the compact set $K$. 
\end{theorem}

From these asymptotics, we can give a simple argument that there
is no limit point of zeros outside the closed unit disk.
Let $\delta>0$ be given. Suppose
$\{x_n\}$ is a sequence of zeros; that is, $F_n(x_n)=0$, that converges
to $x^*$, say, and that $|x_n| \geq 1+\delta$ for all $n$. Then by
Theorem \ref{fnx-out}, 
\[
0= \frac{F_n(x_n)}{ x_n^n}= P(1, 1/x_n) + O( |x_n|^{(C -1) n}).
\]
Since $P(1,1/x^*) \neq 0$, we obtain a contradiction
since $C<1$.
Hence, the zero attractor must lie inside the closed unit disk
$\overline{ {\mathbb D}}$.

 We find that the arguments of the zeros of $F_n(x)$  are
uniformly distributed around the unit circle
by writing $F_n(x)$ as $xg_n(x)$ and applying 
the result of Erd\"os-Tur\'an \cite{erdos} 
(see  equation (\ref{eq:erdos})). Note that 
$g_n(1)=p_n$ and $g_n(x)$ is monic and $g_n(0)=1$.
 We state this result formally as:

\begin{theorem}
The arguments of the zeros of $\{F_n(x)\}$ are uniformly distributed on
the unit circle as $n\to\infty$.
\end{theorem}

Understanding the behavior of zeros inside the unit disk $\mathbb D$
requires a more detailed analysis using the Hardy-Ramanujan circle method.
 A difficulty to
overcome is that the functional equation of the modular
function is unavailable for the generating function
$P(x,u)$.
An important first step in applying the circle method is to rewrite the generating
function $P(x,u)$, for $|x|<1$ fixed,  in a neighborhood of a rational point
$e^{2 \pi i h/k}$
  inside the unit disk $\mathbb D$
  where $h$ and $k$
are relatively prime integers.  Write $u$ as 
$e^{ 2 \pi i (h/k + i z)}$ with $\Re(z)>0$  small. 
The factorization below required careful estimates with
$L$-functions:
\[
\ln [P(x, e^{ 2 \pi i (h/k + i z)})]
=
e^{ w_{h,k}} e^{ \Psi(z)} e^{ j_{h,k}},
\]
where
\begin{eqnarray*}
w_{h,k}
&=&
\frac{1}{2k}\ln(1-x^k)+
\sum_{ \ell, \ell \nmid k} 
\frac{x^\ell}{\ell} \, 
\frac{1}{ e^{- 2 \pi i \ell h/k}-1}, \quad (h,k)=1,
\\
\Psi(z)
&=&
\frac{ \Li_2(x^k)}{ 2 \pi k^2} \, \frac{1}{z},
\\
j_{h,k}(z)
&=&
\frac{1}{ 2 \pi i}
\int_{ -3/4 - i \infty}^{ -3/4 + i \infty} Q_{h,k}(s) \Gamma(s) (2 \pi z)^{-s} \, ds,
\end{eqnarray*}
where 
 $\Li_2(x)$ is the dilogarithm function given on $\mathbb D$ as
$\sum_{n=1}^\infty x^2/n^2$ (see \cite{andrews_sp_fct}, p. 102).
and
$Q_{h,k}(s)$ is defined by means of a series expansion for
$\Re(s)\geq \sigma_0>1$:
\[
Q_{h,k}(s)= \sum_{m \geq 1} \sum_{\ell \geq 1}
\frac{ x^\ell}{\ell} e^{ 2 \pi i \ell m h/k} (\ell m)^{-s}
\]
which
admits an analytic continuation to $\mathbb C$ with a unique singularity,
a simple pole
at $z=1$.

Next we  introduce the quantities needed for the asymptotic
expansion for the polynomials $F_n(x)$:
\[
I_k = 
\frac{1}{  \sqrt{\pi}} \, \frac{1}{ n^{3/4}} \,
\left[ \frac{ \sqrt{ \Li_2(x^k)}  }  {k} \right]^{1/2} \,
\exp
\left(  2 \sqrt{n}  \frac{ \sqrt{ \Li_2(x^k)} }  {k} \right).
\]

\begin{theorem}
Let $K$ be a compact subset of the open upper unit disk. Then
the partition polynomials $F_n(x)$ have the asymptotic form
\begin{eqnarray*}
\lefteqn
{
F_{n}(x)=
e^{w_{0,1}} I_{1} + (-1)^n e^{w_{1,2}} I_{1,2} + e^{-2\pi i n/3} e^{ w_{1,3}} I_{3}
}\\
&&
\qquad \quad + e^{-4 \pi i n/3} e^{w_{2,3}} I_{3}
+
o\left(I_{1}+ I_{2}+I_{3} \right).
\end{eqnarray*}
uniformly on $K$.
\end{theorem}

For simplicity, it is enough to give the zero attractor in the upper unit
disk ${\mathbb D}^+$ since the coefficients of $F_n(x)$ are all real.
Introduce the non-negative subharmonic functions 
$f_k(x) = \frac{1}{k} \Re[ \sqrt{ \Li_2(x^k)}]$
for $|x| \leq 1$.
Let $\mathcal{R}(k)$ be the subset of 
${\mathbb D}^+$ given by:
\[
\mathcal{R}(k) = \{ z  \in {\mathbb D}^+ : |z|  \leq 1, f_k(z) > f_j(z), j \in \{1,2,3\} , j \neq k\}.
\]
Using the same argument as above, we can easily show that there are 
no limit point of zeros that lies in any of the three regions ${\mathcal R}(1)$,
${\mathcal R}(2)$, or ${\mathcal R}(3)$. In fact, the zero attractor consists of the
boundaries of these regions. 
See Figure 1a.
To describe them,  let $C_{k,\ell}$ be the curves
given by $f_k(x) = f_\ell(x)$, $x \in {\mathbb D}^+$ and their subcurves
$\gamma_{k,\ell} = 
C_{k,\ell} \cap \overline{{\mathcal R}(k)}$, where $1 \leq k < \ell \leq 3$. 

\begin{theorem}
The zero attractor of $F_n(x)$ in the upper half-plane
consists of the unit semi-circle  together with the three curves
$\gamma_{k,\ell}$, $1 \leq k < \ell \leq 3$.
\end{theorem}

\begin{figure}\label{fig:stanley}
\begin{center}
\includegraphics[height=4.5cm,width=4.5cm]{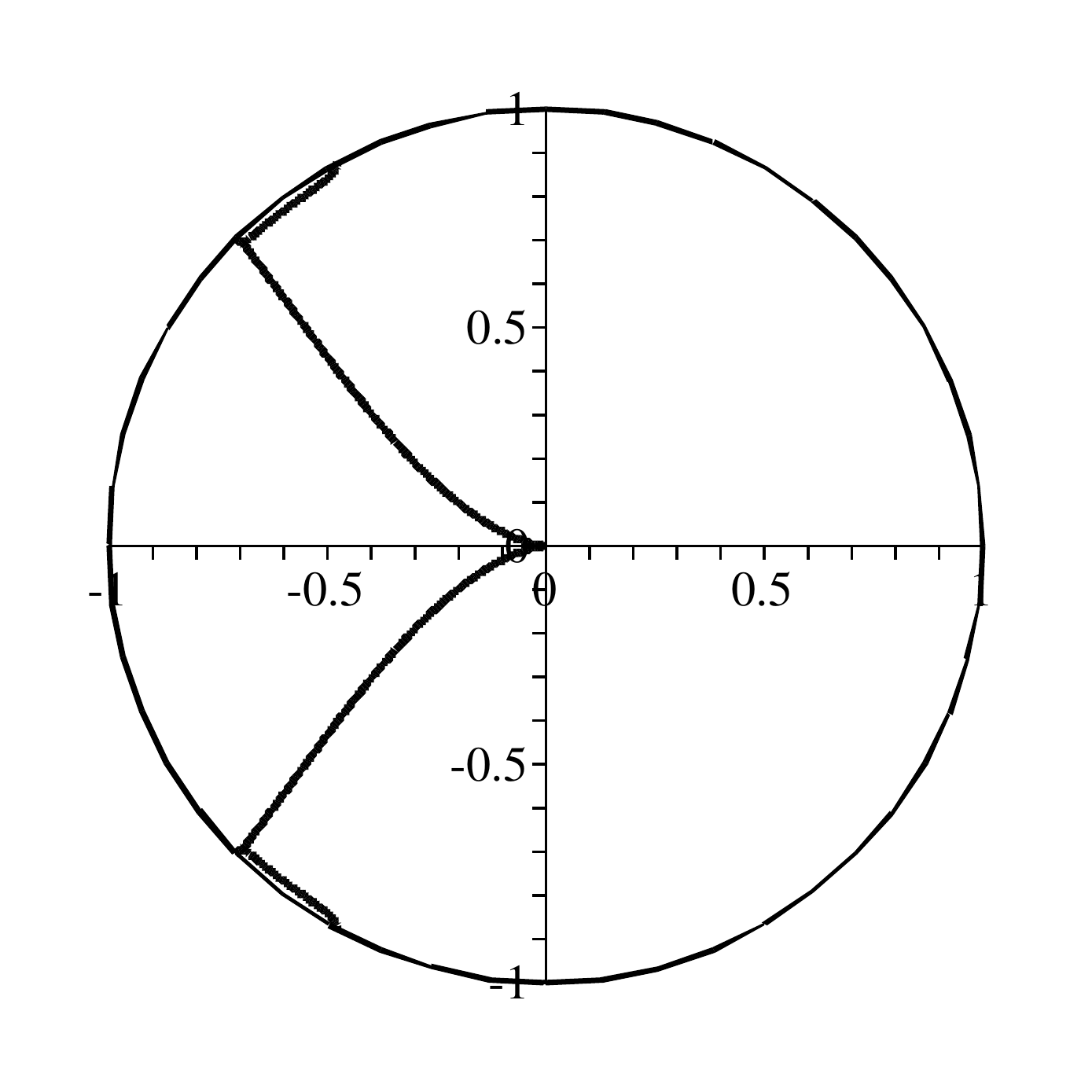} 
\quad
\includegraphics[height=4.5cm,width=4.5cm]{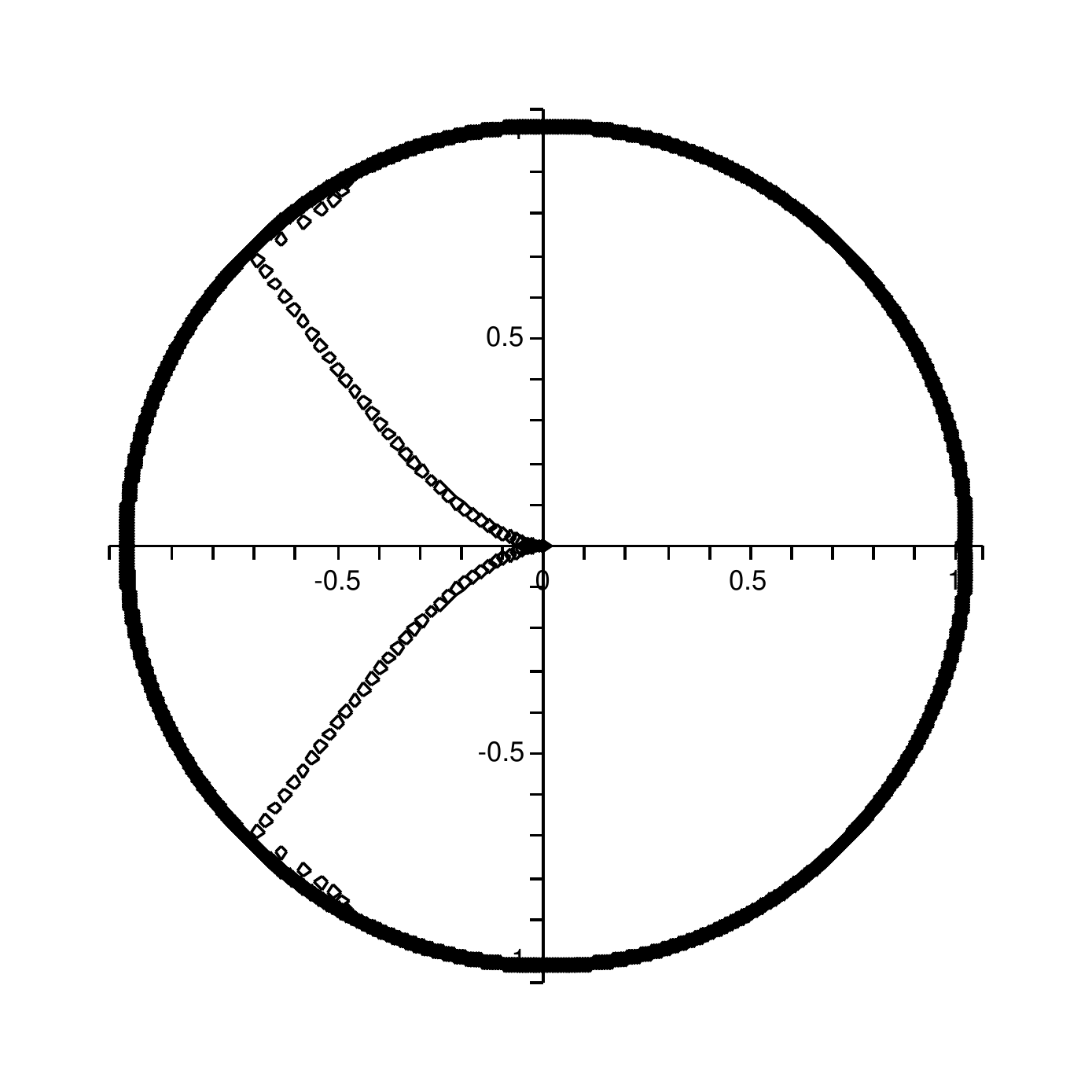} 
\caption{(a) Zero Attractor for Partitions with Parts Polynomial;
(b) All Zeros of $F_n(x)$, $n=10,000$}
\end{center}\end{figure}

A basic estimation of the number of zeros of $F_n(x)$ relative to the
unit disk $\mathbb D$ demonstrates a striking dichotomy.

\begin{theorem}
(a)
 Let $\varepsilon >0$. Then
$
 \# \{ z : F_n(x)=0, 1 \leq |z| \leq 1 + \varepsilon \} = O(n).
$
\\
(b)
Let $K$ be any compact subset of the open unit disk $\mathbb D$. Then
\[
\#  (Z(F_n) \cap K) = O(\sqrt{n}).
\]
\end{theorem}

\begin{corollary}
The asymptotic zero distribution for $\{F_n(x)\}$ is Lebesgue measure on the
unit circle.
\end{corollary}

Clearly, the standard definition for the asymptotic zero distribution 
(Definition \ref{def:asymptotic_zero})
ignores any contribution from the zeros inside the unit disk $\mathbb D$.
As a consequence, it is necessary to extend that definition:

\begin{definition}
The {\sl asymptotic zero distribution of order $\alpha$ on
a domain $D$}
for a sequence $\{q_n\}$ of polynomials 
whose degrees go to $\infty$
is the weak$^{*}$-limit $\mu$
of the normalized counting measures of their zeros
$\frac{1}{\deg(q_n)^\alpha} \sum \{ \delta_z : z \in {\mathcal Z}(q_n)
\textrm{  and  } z \in D\}$.
\end{definition}

For the sake of exposition, we will restrict our discussion of the
order $1/2$ asymptotic zero distribution
 $\mu$ for $\{F_n(x)\}$ to the upper unit disk
${\mathbb D}^+$. Since the zero attractor consists of three analytic
curves, 
the support of the measure $\mu$ is supported exactly on those curves;
in particular,
it  will be enough to describe $\mu$ in a neighborhood of each
of them.

\begin{theorem}
Let $1\leq k<\ell \leq 3$. For each curve
$\gamma_{k,\ell}$
in the zero attractor
 there exists a neighborhood $U_{k,\ell}$ of $\gamma_{k,\ell}$
 and a conformal map $G_{k,\ell}$ on $U_{k,\ell}$ that maps 
 $\gamma_{k,\ell}$ into the unit circle such that the 
 asymptotic zero distribution of $\{F_n(x)\}$ in $U_{k,\ell}$
 of order $1/2$ is the pull-back of Lebesgue measure on the unit circle.
\end{theorem}

For the specifics of these mappings, see \cite{boyer_goh}. 
Their construction comes from the explicit asymptotic expansion of the polynomials $F_n(x)$.

\section{Rank and Crank Polynomials}\label{section:crank}

We now emphasize another way to look at the partition in parts polynomials $F_n(x)$
relative to their generating function
to show its similarities with other generating functions that appear in partition theory.
 It is well known (see \cite{andrews_sp_fct}, p. 568)
that
\[
P(z,q)
=
1+\sum_{n=1}^\infty F_n(z) q^n
=
1+ \sum_{n=1}^\infty \frac{  q^{n^2} z^n }{  (q;q)_n (zq;q)_n}
= \frac{1}{(zq;q)_\infty};
\]
namely, $F_n(x)$ is the coefficient of $q^n$ of the generating function.
With this viewpoint, there are several other natural polynomial families 
defined in terms of 
either Durfee squares (see \cite{durfee_poly}),
ranks, or cranks. 
Here we are using the standard notations
$(q)_n=(1-q)(1-q^2)\cdots (1-q^{n-1})$ and, more generally,
$(a;q)_n = (1-a)(1-aq) \cdots (1-aq^{n-1})$;
next, when $|q|<1$, let
 $(a;q)_\infty = \lim_{n \to \infty} (a;q)_n$ and
 $(q)_\infty$ for $(1;q)_\infty$.

For a partition $\lambda$ of $n$, its Durfee square is the largest 
square that lies inside its Ferrers graph (see \cite{andrews},
Chapter 2).
The polynomials $d_n(z)$ for Durfee squares were introduced
in \cite{durfee_poly} and are given in
terms of their generating function:
\[
D(z,q)= 
\sum_{n=1,k=1}^\infty d(n,k) q^n z^k
=
\sum_{n=1}^\infty d_n(z) q^n
=
 \sum_{n=1}^\infty \frac{ q^{n^2} z^n}{ (q;q)_n^2} 
\]
where $d(n,k)$ is the number of partitions of $n$ with a Durfee square of size $k$. Further, in \cite{durfee_poly} and \cite{canfield}, they conjecture that the associated polynomials $\{d_n(z)\}$ have only negative real zeros. Note that the Erd\"os-Tur\'an result does not apply here since the degree of $d_n(z)$ is $\lfloor\sqrt{n} \rfloor$.

F. Dyson introduced the statistic of rank for a partition $\lambda$
of $n$ as the difference between its largest part and the number
of its parts (see \cite{andrews}, p. 142).
We introduce the  rank polynomials $r_n(z)$  as follows.
Consider their generating function: 
\[
R(z,q)= \sum_{n=0}^\infty \, \sum_{m=-\infty}^\infty
N(m,n) z^m q^n 
=
\sum_{n=0}^\infty \, \sum_{m=-(n-1)}^{n-1}
N_n(z) q^n
=
1+ \sum_{n=1}^\infty \frac{ q^{n^2}}{ (zq;q)_n (z^{-1}q;q)_n}
\]
where $N(m,n)$ is the number of partitions of $n$ with rank $m$
and $\{N_n(z)\}$ are symmetric Laurent polynomials.
Set $r_n(z) $ to be the principal part of $N_n(z)$
and call it the {\sl rank polynomial}:
\[
r_n(z) = \sum_{m=0}^{n-1} N(m,n) z^m.
\]

Let  $\lambda$ be the partition of $n$  given as
$ \lambda_1 + \cdots + \lambda_s
+ 1 + \cdots 1$, where there are exactly $r$ 1's.
Let $o(\lambda)$ be the number of parts $>r$.
Then the crank of $\lambda$ is $\lambda_1$ if $r=0$
and $o(\lambda)-r$ if $r > 0$  \cite{andrews_garvan}.
Let $M(m,n)$ be the number of partitions of $n$ whose crank
is exactly $m$. Then $C(z,q)$ is their generating function
for $n>1$ where
\begin{eqnarray*}
\lefteqn
{
C(z,q)= \sum_{n=0}^\infty \,  \sum_{m=-\infty}^\infty
M(m,n) z^m q^n
=
\sum_{n=0}^\infty \,  \sum_{m=-\infty}^\infty
M_n(z) q^n
}\\
&&
\qquad 
=
\prod_{n=1}^\infty \frac{ (1-q^n)}{ (1-zq^n) (1-z^{-1}q^n)}
=
\frac{ (q)_\infty}{ (z;q)_\infty (z^{-1};q)_\infty}.
\end{eqnarray*}
Let  $c_n(z)$ be the principal part of $M_n(z)$
and call it the
 {\sl crank polynomial}:
\[
c_n(z) = \sum_{m=0}^{n} M(m,n) z^m.
\]

We can apply the Erd\"os-Tur\'an result on the asymptotic
distribution of the arguments of the zeros to the two families
for the rank and crank polynomials
since their coefficients are all non-negative, they are monic,
and both quotients
$r_n(1)/\sqrt{r_n(0)}$, $c_n(1)/\sqrt{c_n(0)}$
are bounded above by  $\sqrt{p_n}$.
We record this as a theorem:

\begin{theorem}
The arguments of the zeros of both the rank and crank polynomials
are  uniformly distributed on the unit circle as $n\to \infty$.
\end{theorem}

\begin{figure}\label{fig:rank_zeros}
\begin{center}
\includegraphics[height=4.5cm,width=4.5cm]{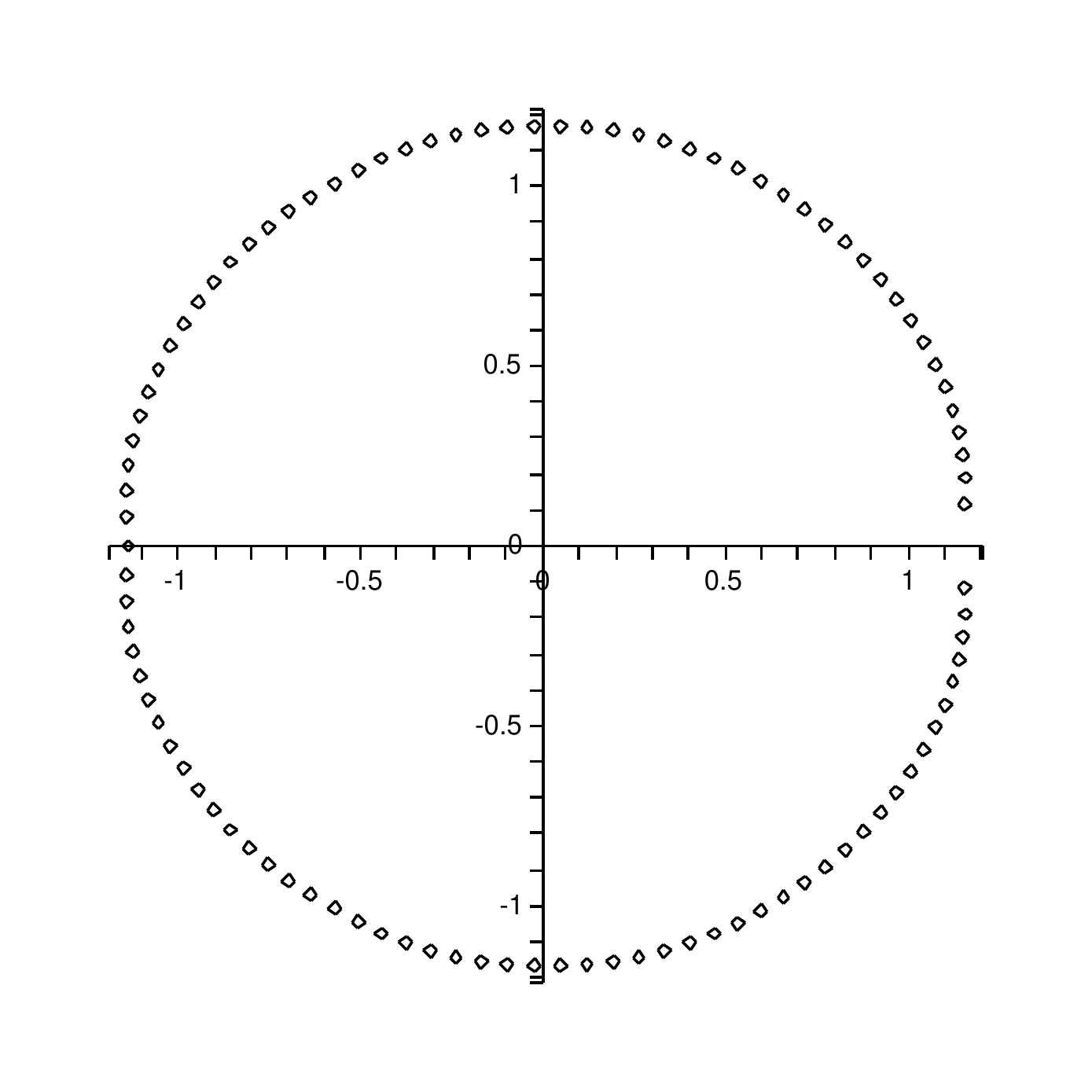} 
\includegraphics[height=4.5cm,width=4.5cm]{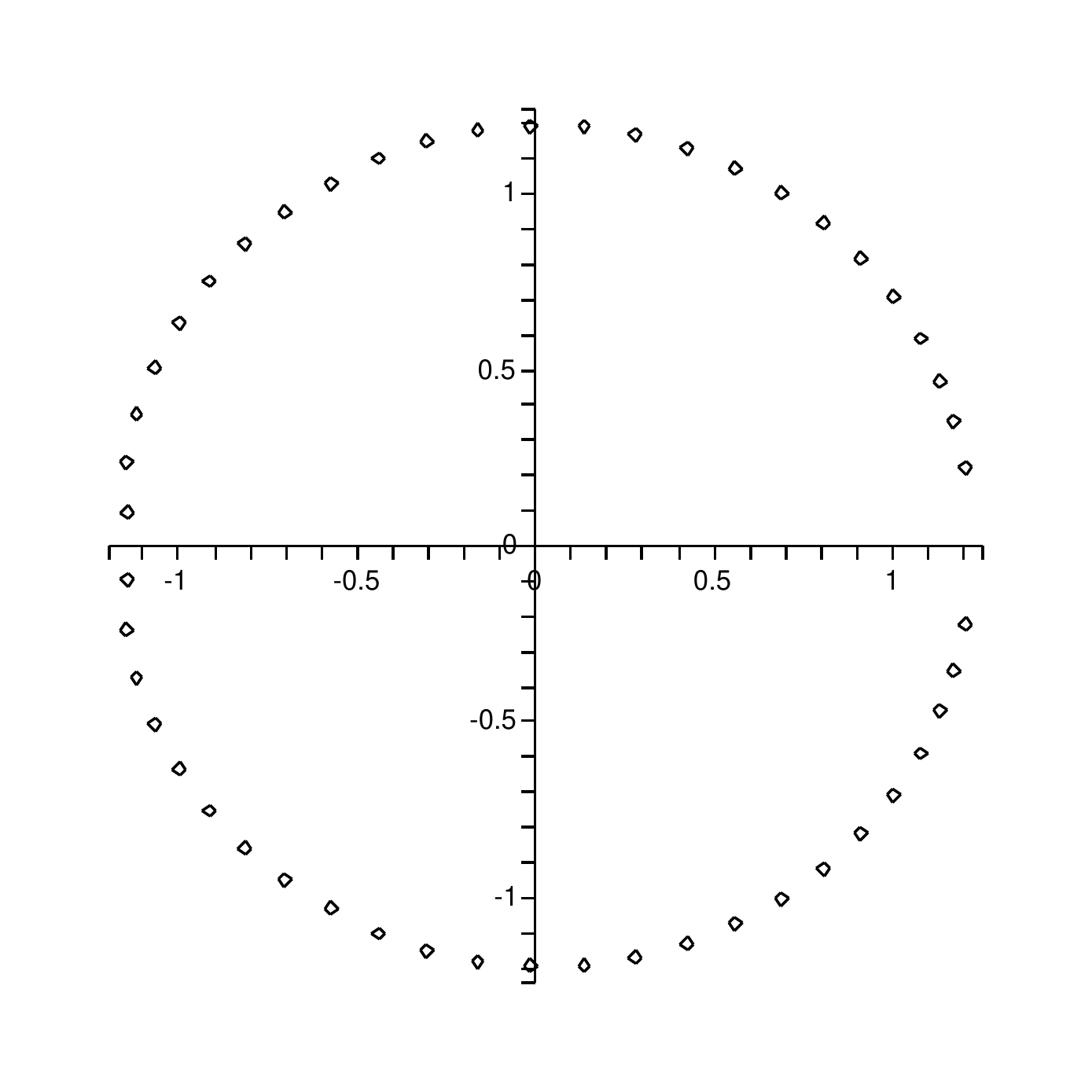} 
\caption{Zeros of the rank  polynomial degree 100;
crank polynomial degree 50}
\end{center}
\end{figure}

From explicit computation, we find that their zero attractor appears to
be the unit circle (see Figure 2).
It is very natural  to attempt to extend the work   for the partition in parts
polynomials $F_n(x)$ to establish this conjecture.

Furthermore, 
it would be interesting to see how any  partition
polynomials in this paper  fit into the statistical mechanics framework described in
Vershik's paper \cite{vershik}.

\section{Summary}

For all but one of the partition polynomial families, 
the unit circle has a dominant role.
Their zero attractor is either equal or contains the unit circle while their
asymptotic zero distribution involves Lebesgue measure on the unit circle.
All this makes it even more intriguing to understand the meaning of
the subtle two-scale asymptotics of the partition in parts polynomials
$F_n(x)$ in Section \ref{section:stanley}.

\vspace{0.3in}
\noindent
Department of Mathematics
\\
Drexel University
\\
Philadelphia, PA 19104
\\
{\sl email}: {\tt rboyer  at   math.drexel.edu}

\end{document}